# HOMOCLINIC ORBITS IN A POROUS PELLET

Andrzej Burghardt, Marek Berezowski
Institute of Chemical Engineering, Polish Academy of Sciences, Gliwice, Poland
E-mail: marek.berezowski@polsl.pl

**Abstract**

The analysis performed as well as extensive numerical simulations have revealed the possibility of the generation of homoclinic orbits as a result of homoclinic bifurcation in a porous pellet. A method has been proposed for the development of a special type of diagrams - the so-called bifurcation diagrams. These diagrams comprise the locus of homoclinic orbits together with the lines of limit points bounding the region of multiple steady states as well as the locus of the points of Hopf bifurcation. Thus, they define a set of parameters for which homoclinic bifurcation can take place. They also make it possible to determine conditions under which homoclinic orbits are generated.

Two kinds of homoclinic orbits have been observed, namely semistable and unstable orbits. It is found that the character of the homoclinic orbit depends on the stability features of the limit cycle which is linked with the saddle point.

Very interesting dynamic phenomena are associated with the two kinds of homoclinic orbits; these phenomena have been illustrated in the solution diagrams and phase diagrams.

## 1. INTRODUCTION

The coupling of physical and chemical processes such as diffusional and convective mass and heat transport and chemical reactions may lead, in a given system, to a number of interesting dynamic phenomena which can be rather difficult to predict without a thorough mathematical analysis. These phenomena include the multiplicity of steady states and the accompanying hysteresis, as well as the oscillatory destabilization of the system, which can also lead to the appearance of hysteresis.

All these phenomena are of primary importance in the design, start-up and control of chemical reactors. If these phenomena are insufficiently described either qualitatively (no information about the possible existence of multiple steady states) or quantitatively (unknown amplitude and frequency of self-excited oscillations of the system), the reactor may not operate at its optimum or can even be severely damaged (Vleeshouver et al., 1992, Morud et al., 1998, Mancusi et al., 2000).

In heterogeneous catalytic reactors the chemical transformation of matter takes place in a catalyst pellet. Consequently, all the processes and phenomena associated with this transformation are generated within the pellet. The fluid flowing through the reactor supplies reactants for the reactions occurring in the pellet and takes away the products formed therein. Obviously, the fluid affects the processes within the pellet, but only in an indirect way, through the appropriate boundary conditions at the fluid-pellet interface. Thus, the process occurring in a heterogeneous catalytic reactor is determined by the phenomena within the catalyst pellet. A precise quantitative description of these phenomena is therefore of utmost importance.

The problem of defining the exact criteria determining the conditions of the occurrence of steady-state multiplicity in porous catalyst pellets has received a great deal of attention and may now be regarded as solved (Hu, Balakotaiah and Luss, 1985a, 1985b, 1986, Burghardt and Berezowski, 1989, 1990a, 1990b).

Consequently, the question arises about the stability of the stationary states thus determined, especially in the multiplicity region. Analysis of the stability of steady states of a catalyst pellet was performed by means of the various approximation methods such as the linearization technique (Hlavacek et al., 1969a, 1969b) or the averaging technique (Luss and Lee, 1971, Lee et al., 1972 and Ray and Hastings, 1980).

Burghardt and Berezowski (1991) analysed the stability of stationary solutions of the system considered with respect to static bifurcations, i.e., the analysis concerned the case of small Lewis numbers (Le → 0). As a result of partial integration and further transformation, the system of two partial differential equations has been reduced to one ordinary differential equation. The stability analysis of this equation led to the characteristics of the static bifurcation by means of developing analytical equations for the local and global stability.

A number of authors (Sony, Schmidt and Aris, 1991, Van den Bosch and Padmanabhan, 1974, Subramanian and Balakotaiah, 1996) have placed emphasis on deriving a method whereby the parameter region could be determined in which the bifurcation to periodic solutions is possible. The most important result of the studies by Burghardt and Berezowski (1993) is a set of analytical relations between the parameters of the system which define the limits of oscillatory instability (Hopf bifurcation) and saddle type instability as well as approximate formulae to determine the critical values of the Lewis number (i.e. the value above which an oscillatory destabilization of the system becomes possible).

The possibility of determining the limits of oscillatory instability without resorting to tedious and time-consuming computations is especially useful in modelling processes that occur in heterogeneous chemical reactors.

Of special interest are the conclusions of the study by Burghardt and Berezowski (1997), where extensive numerical results have revealed three scenarios for the generation of self-excited oscillations of the temperature and composition inside the pellet. Besides the well-known supercritical and subcritical Hopf bifurcations to oscillatory solutions the authors also observed homoclinic bifurcation in the region of multiple steady states. To the best of our knowledge the possibility of the homoclinic bifurcation occurring in the dynamic processes taking place in a porous catalyst pellet has not so far been reported in the existing literature.

The generation of a homoclinic orbit proceeds according to the following scenario. Changing the bifurcation parameter, part of an existing limit cycle moves closer and closer to a saddle point. At the homoclinic bifurcation the limit cycle touches the saddle point and becomes a homoclinic orbit. According to Strogatz (1994) and Wiggins (1988) this phenomenon represents another type of the infinite-period or global bifurcation. The key feature of this bifurcation is the fact that the stable and unstable manifolds of the saddle point join the existing limit cycle, and thus the homoclinic orbit connects the saddle point to itself, developing a limit cycle of an infinite period (Guckenheimer and Holmes, 1983, Thompson and Stewart, 1993, Wiggins, 1998).

The aim of the analysis carried out in the present study is to determine conditions under which homoclinic bifurcation can develop, using the model of the processes occurring in the catalyst pellet. The next step is to investigate the dynamics of the formation of the homoclinic orbit.

It has to be borne in mind that these goals cannot be achieved based solely on the analysis of the dynamics of the generation of the homoclinic orbit following the change in a selected bifurcation parameter of the model. Rather, this analysis should be carried out in relation with other bifurcation phenomena, such as bifurcation to multiple steady state solutions and Hopf bifurcation.

## 2. MODEL OF THE PROCESS AND A METHOD FOR DEVELOPING BIFURCATION DIAGRAMS AND PHASE DIAGRAMS

The basic equations describing the non-stationary process of mass and heat transport with a simultaneous chemical reaction in a porous catalyst pellet are as follows:

$$\frac{\partial C_A}{\partial t} = D_e \frac{1}{\rho^a} \frac{\partial}{\partial \rho}\left(\rho^a \frac{\partial C_A}{\partial \rho}\right) + \nu_A r(C_A, T) \tag{1}$$

$$(\rho c)_k \frac{\partial T}{\partial t} = \lambda_e \frac{1}{\rho^a} \frac{\partial}{\partial \rho}\left(\rho^a \frac{\partial T}{\partial \rho}\right) + (-\Delta H) r(C_A, T) \tag{2}$$

The transport across the boundary is determined by the following relationships between the fluxes:

$$\text{at} \quad \rho = \rho_0 \qquad k(C_A - C_{A0}) = -D_e \frac{\partial C_A}{\partial \rho} \tag{3}$$

$$\alpha(T - T_0) = -\lambda_e \frac{\partial T}{\partial \rho} \tag{4}$$

If we restrict ourselves to symmetrical profiles, the remaining boundary conditions take the form:

$$\text{at} \quad \rho = 0 \qquad \frac{\partial C_A}{\partial \rho} = 0 \quad \text{and} \quad \frac{\partial T}{\partial \rho} = 0 \tag{5}$$

The transient eqs. (1) and (2) as well as boundary conditions (3)-(5) can be expressed in the following dimensionless form:

$$\frac{\partial y}{\partial \tau} = \frac{1}{x^a} \frac{\partial}{\partial x}\left(x^a \frac{\partial y}{\partial x}\right) + \nu_A \phi_0^2 R(y, z) \tag{6}$$

$$\frac{1}{Le'} \frac{\partial z}{\partial \tau} = \frac{1}{x^a} \frac{\partial}{\partial x}\left(x^a \frac{\partial z}{\partial x}\right) + \beta \phi_0^2 R(y, z) \tag{7}$$

subject to boundary conditions:

$$\text{at} \quad x = 0 \qquad \frac{\partial y}{\partial x} = 0 \quad \text{and} \quad \frac{\partial z}{\partial x} = 0 \tag{8}$$

$$\text{at} \quad x = 1 \qquad y = 1 - \left(\frac{1}{Bi_M}\right)\left(\frac{\partial y}{\partial x}\right) \tag{9}$$

$$\text{and} \quad z = 1 - \left(\frac{1}{Bi_H}\right)\left(\frac{\partial z}{\partial x}\right) \tag{10}$$

The dimensionless variables together with the dimensionless numbers are defined in "Notation".

In our previous analyses concerning the stability of steady-state solutions (Burghardt and Berezowski, 1993, 1997) we employed a model of a porous catalyst pellet which takes into account the concentration gradients inside the pellet and assumes that the pellet temperature is uniform but different from that of the surrounding fluid. Thus, the resistance to heat transfer is concentrated in a boundary layer around the pellet, while the resistance to mass transfer occurs only inside the pellet.

The assumption concerning the mass-transfer resistance $Bi_M \to \infty$ transforms boundary condition (9) to the following form:

at $\quad x = 1 \qquad y \approx 1$ $\hfill (11)$

Before introducing the second assumption (concerning the heat transfer), we average Eq. (7) over the volume of the pellet:

$$\frac{1}{Le'}\frac{1}{V}\int_V \frac{\partial z}{\partial \tau}dV = \frac{1}{V}\int_V \frac{1}{x^a}\frac{\partial}{\partial x}\left(x^a \frac{\partial z}{\partial x}\right)dV + \beta\phi_0^2 \frac{1}{V}\int_V R(y,z)dV \tag{12}$$

Introducing the following definitions:

$$\frac{1}{V}\int_V R(y,z)dV = \overline{R}(y,z) \tag{13}$$

and

$$\frac{1}{V}\int_V \frac{\partial z}{\partial \tau}dV = \frac{\partial}{\partial \tau}\frac{1}{V}\int_V z\,dV = \frac{\partial \overline{z}}{\partial \tau} \tag{14}$$

and expressing the value of $dV/V$ for the three pellet shapes (slab, cylinder, sphere) as

$$\frac{dV}{V} = (a+1)x^a dx \tag{15}$$

we can rewrite Eq. (12) in the following form:

$$\frac{1}{Le'}\frac{\partial \overline{z}}{\partial \tau} = (a+1)\left(\frac{\partial z}{\partial x}\bigg|_{x=1} - \frac{\partial z}{\partial x}\bigg|_{x=0}\right) + \beta\phi_0^2 \overline{R}(y,z) \tag{16}$$

Using boundary conditions (8) and (10) we finally obtain an equation defining the mean temperature of the catalyst pellet:

$$\frac{1}{Le'}\frac{\partial \overline{z}}{\partial \tau} = (a+1)Bi_H(1-z)_{x=1} + \beta\phi_0^2 \overline{R}(y,z) \tag{17}$$

Introducing the second assumption, according to which the intraparticle temperature gradient is negligible in comparison with the temperature difference in the fluid film surrounding the pellet, we can write

$$\overline{z} = z\big|_{x=1} = z(\tau) \tag{18}$$

Finally, the system of equations describing the process according to the model assumed takes the form:

$$\frac{\partial y}{\partial \tau} = \frac{1}{x^a}\frac{\partial}{\partial x}\left(x^a \frac{\partial y}{\partial x}\right) + v_A \phi_0^2 R(y,z) \tag{19}$$

$$\frac{1}{Le'}\frac{\partial z}{\partial \tau} = (a+1)Bi_H(1-z) + \beta\phi_0^2 \overline{R}(y,z) \tag{20}$$

subject to boundary conditions:

at $\quad x = 1 \qquad y = 1$ $\hfill (21)$

at $\quad x = 0 \qquad \frac{\partial y}{\partial x} = 0$ $\hfill (22)$

and initial conditions:

$\tau = 0 \qquad y(x) = y_0(x)$ $\hfill (23)$

$\tau = 0 \qquad z = z_0$ $\hfill (24)$

Introducing, as previously, the modified Thiele modulus valid for power law kinetics

$$\theta_0 = \frac{\phi_0}{a+1}\left(\frac{n+1}{2}\right)^{1/2} \tag{25}$$

and transforming the energy balance equation (20), the following set of equations can be obtained:

$$\frac{\partial y}{\partial \tau} = \frac{1}{x^a}\frac{\partial}{\partial x}\left(x^a \frac{\partial y}{\partial x}\right) + q(a,n)\theta_0^2 R(y,z) \tag{26}$$

$$\frac{1}{Le}\frac{dz}{d\tau} = 1 - z + \beta^*\theta_0^2\overline{R}(y,z) \tag{27}$$

subject to boundary conditions (21) and (22) and initial conditions (23) and (24). The modified parameters Le and $\beta^*$ are defined as follows:

$$Le = Le'Bi_H(a+1) = \frac{\alpha L(a+1)}{(\rho c)_k D_e} \tag{28}$$

and

$$\beta^* = \frac{2\beta(a+1)}{Bi_H(n+1)} = \frac{(-\Delta H)C_{A0}D_e}{\alpha L T_0}\frac{2(a+1)}{n+1} \tag{29}$$

where

$$q(a,n) = \frac{2(a+1)^2}{n+1}\nu_A \tag{30}$$

Assuming the following relation describing the reaction rate

$$R(y,z) = \exp\left[\gamma\left(1 - \frac{1}{z}\right)\right]y^n \tag{31}$$

we replace z with a new state variable

$$\theta^2 = \theta_0^2 \exp\left[\gamma\left(1 - \frac{1}{z}\right)\right] \tag{32}$$

which is more convenient in further analysis.

For the appropriate representation of the formation of homoclinic orbits and the region of their existence it is of crucial importance to illustrate, in a single bifurcation diagram (Kubiĉek and Marek, 1983), the coexistence of three basic bifurcation phenomena occurring in the catalyst pellet. In the present study these coexisting phenomena are defined as bifurcation to multiple steady state solutions for the pellet (LP), Hopf bifurcation (HB) and homoclinic bifurcation (Hcl). Also, solution diagrams (Kubiĉek and Marek, 1983) are derived which show the branches of both steady and unsteady states of the pellet, the points of the Hopf bifurcation and the branches of periodic solutions along with their complete amplitudes. These diagrams illustrate, for a strictly defined set of the model parameters, the process of generating the homoclinic bifurcation. The detailed location of the homoclinic orbit relative to the points determining steady states of the pellet and the possible limit cycles is shown in the phase diagram in the system of coordinates which define state variables of the system analysed as functions of time.

The approach used to derive these three types of diagrams is as follows.

Using the method of orthogonal collocation the model equations for the pellet are discretized, thus leading to n ordinary differential equations with temporal derivatives. This system is then analysed numerically.

Based on the one-parameter continuation procedure solution diagrams of the steady states are derived along with the points of the Hopf bifurcation (HB) situated on their branches. These diagrams form the basis for the construction of bifurcation diagrams. The bifurcation diagrams illustrate the regions of multiple steady states bounded by the lines of static bifurcations (LP), as well as the regions of oscillations bounded by the lines of the Hopf bifurcation (HB). The two sets of bounding lines are obtained employing the two-parameter continuation procedure.

The calculations were performed using the Auto 97 package (Doedel et al., 1997).

In order to obtain a full picture of the branches of periodic solutions (i.e. maximums and minimums of the amplitudes) necessary in further analysis, the package had to be modified (originally, it generated only maximums of these solutions).

The most complex task was to determine numerically the lines of homoclinic bifurcations (Hcl) in the bifurcation diagram. Therefore, the following procedure was adopted. The solution diagram was determined in which the branch of periodic solutions was linked with the saddle point on an unstable branch of steady states thus forming the homoclinic orbit. Parameters of a single point of the homoclinic bifurcation were thus obtained. Then, starting at this point of the homoclinic bifurcation the line of the homoclinic bifurcation was determined (Hcl) in the bifurcation diagram by using the two-parameter continuation procedure. Of great help were suggestions made by one of the authors of Auto 97, Prof. Eusebius Doedel.

In order to develop the homoclinic orbit in the phase diagram it was necessary to analyse the full dynamics of its generation by means of numerical integration of the system of nonlinear equations (26) and (27) subject to boundary conditions (21) and (22) and different initial conditions (23) and (24).

The integration led to the temperature transients for the catalyst pellet and to the transient concentration profiles inside the pellet. The classification of the dynamic behaviour of a given system may be best done using two-dimensional phase diagrams. Therefore the concentration profile has been replaced with the spatial integral average of this profile:

$$\eta(\tau) = \int_0^1 y(x,\tau) x^a (a+1) dx \qquad (33)$$

If, then, the variables $\eta$ i $\theta$ are used in the phase diagram, it is possible to carry out a rigorous analysis of the dynamic behaviour of the system and to determine the homoclinic orbit which connects the saddle point to itself, as well as the location of other stationary points and possible limit cycles. The exceptionally strong sensitivity of the location of the homoclinic orbit in the phase diagram to the initial conditions used in numerical simulations led to serious problems.

## 3. RESULTS AND DISCUSSION

The bifurcation diagrams discussed earlier were derived using two model parameters as coordinates, namely the parameter $\gamma$ and the Thiele modulus $\theta_0$, at constant values of the parameters $\beta^*$ and Le. In the solution diagrams the dimensionless pellet temperature $\theta$ is selected as a state variable, while the Thiele modulus $\theta_0$ is a bifurcation parameter.

The phase diagrams, which illustrate the dynamic behaviour of the process analysed, are drawn in the system of coordinates $\eta, \theta$ representing state variables of the system.

The analysis was performed for two different values of the Lewis number: Le = 10 and Le = 50, at a constant value of $\beta^* = 0.25$.

In Table 1 the values of the Lewis numbers are listed for the various operating conditions typical of actual processes in fixed-bed reactors. This table fully confirms the practical significance of the choice of the two Lewis numbers employed in our analysis.

Table 1. Lewis numbers for the various operating conditions

| P(MPa) | L = 1 (cm) | | | | | L = 0.5 (cm) | | | | |
|---|---|---|---|---|---|---|---|---|---|---|
| | Re | $\alpha$ $[W/(m^2 K)]$ | $Le = \alpha L/(\rho c)_k \, D_e$ | | | Re | $\alpha$ $[W/(m^2 K)]$ | $Le = \alpha L/(\rho c)_k \, D_e$ | | |
| 0.1 | 401 | 70.2 | 0.60 | 6.00 | 60.0 | 201 | 100.8 | 0.43 | 4.3 | 43 |
| 1.0 | 3973 | 274 | 2.34 | 23.4 | 234 | 1987 | 362.5 | 1.55 | 15.5 | 155 |
| 10.0 | 38,576 | 1105 | 9.44 | 94.4 | 944 | 19,288 | 1286 | 5.92 | 59.2 | 592 |
| 30.0 | 112,453 | 2312 | 19.8 | 198 | 1980 | 56,226 | 3072 | 13.1 | 131 | 1310 |
| | | | $10^{-6}$ | $10^{-7}$ | $10^{-8}$ | | | $10^{-6}$ | $10^{-7}$ | $10^{-8}$ |
| | | | | $D_e \, (m^2/s)$ | | | | | $D_e \, (m^2/s)$ | |

It is seen from Fig. 1, which is the bifurcation diagram for $Le = 10$, that below $\gamma = 7.65$ only single stable steady states exist for any value of $\theta_0$. Over the range $\gamma = 7.650 - 7.939$ the loci of the Hopf bifurcation determine such intervals of the Thiele modulus over which periodic solutions exist. These loci form the limits of the branches of periodic solutions in the solution diagrams, i.e. the points of the Hopf bifurcation (HB).

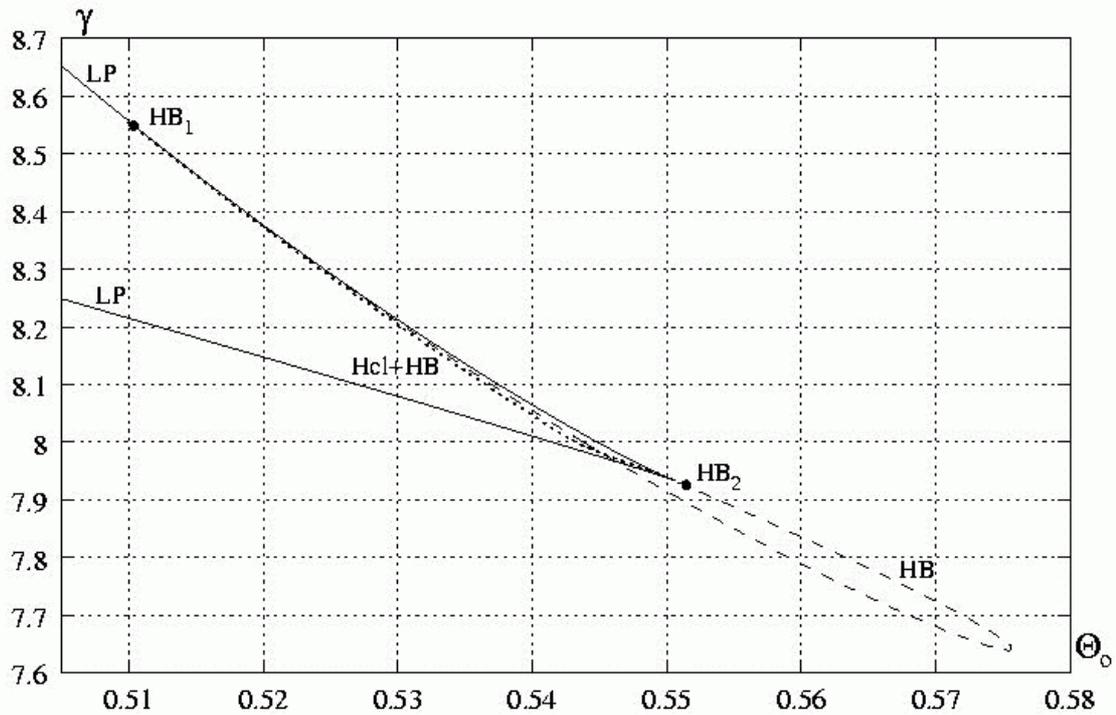

Fig.1. Bifurcation diagram: $a = 0, n = 1, \beta^* = 0.25, Le = 10$. — LP – locus of limit points, – – – HB – locus of Hopf bifurcation points, $\cdots$ Hcl – locus of homoclinic bifurcation points.

An interesting feature from the standpoint of the dynamics of the system studied is the moment at which the locus of the Hopf bifurcations begins to touch the lines limiting the region of multiple solutions (the so-called limit points), and then enters this region. It can be seen from Fig. 2 (which is an enlargement of a part of Fig. 1) that the line of the points of the Hopf bifurcations, once inside the region of multiple solutions, becomes automatically the locus of the homoclinic bifurcation. This is a direct proof that for the values of $\gamma$ and $\theta_0$ which lie on this line the limit cycle has indeed joined the saddle point as a result of the so-called infinite-period bifurcation and the formation of a homoclinic orbit.

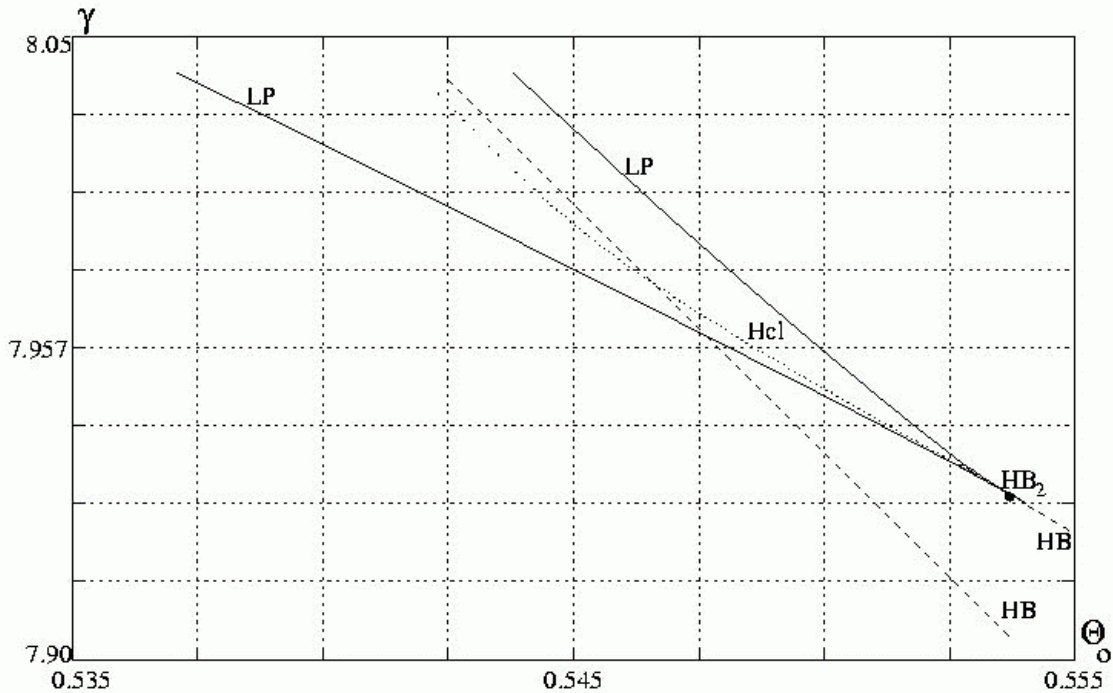

For the values of $\gamma = 7.94 - 7.957$ the branch of periodic solutions extends from the lower line of the points of the Hopf bifurcations up to the line of homoclinic bifurcations inside the region of multiple solutions.

To carry out a detailed analysis of the generation of homoclinic bifurcation and the formation of homoclinic orbit, solution diagrams were prepared (Figs 3 and 4). In Fig. 3 $\left(\gamma = 7.939\right)$ two points of the Hopf bifurcation still exist which constitute the boundary of the branches of periodic solutions; the region of multiple solutions is only being formed, as demonstrated by two barely visible limit points (LP). However, already for $\gamma = 7.940$ the limit cycle clearly meets the saddle point of the unstable branch of steady states and the homoclinic orbit is formed.

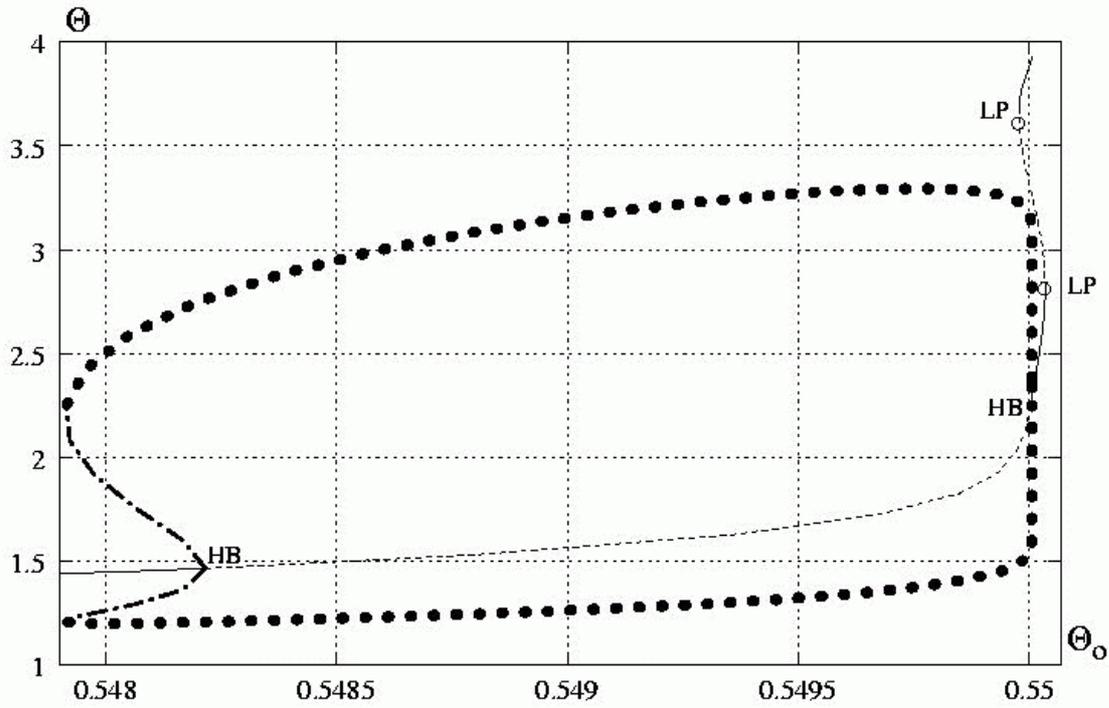

Fig.3.    Solution   diagram:   $a = 0, n = 1, \beta^* = 0{,}25, \gamma = 7{,}939, \mathrm{Le} = 10$ .  — stable   steady   states, – – –
unstable steady states, • • • stable limit cycles, – • – unstable limit cycles.

 Consequently, over the range $\gamma = 7.94 - 7.957$ solution diagrams exist of the form  presented  in
Fig. 4.

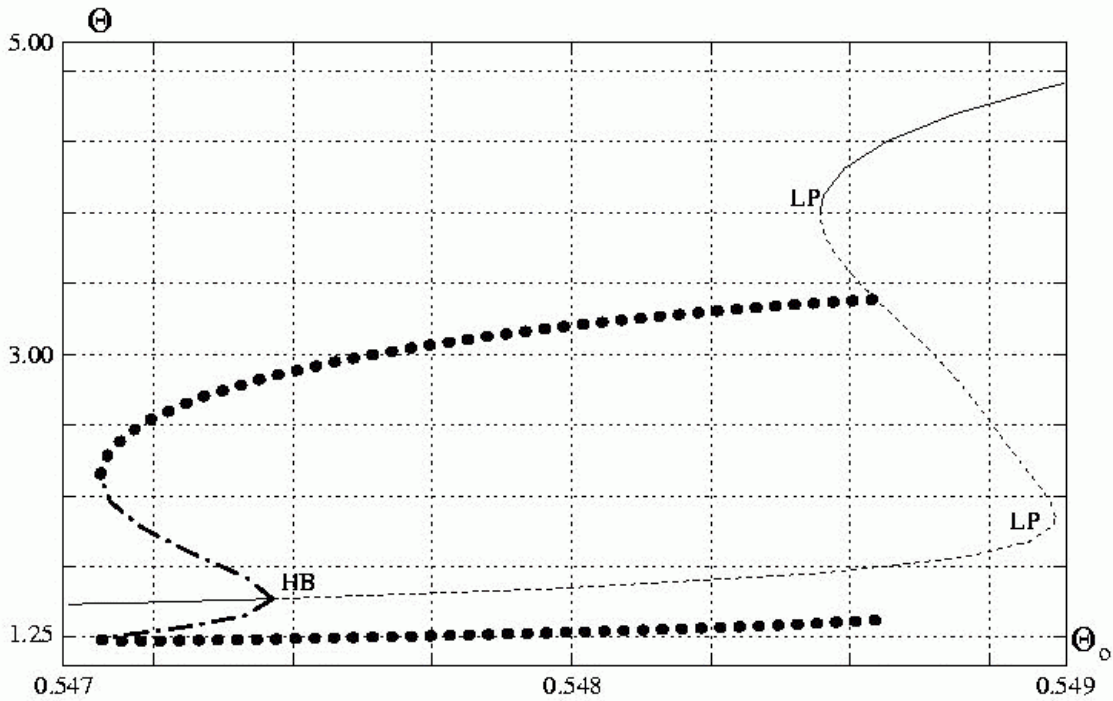

Fig.4. Solution diagram: $a = 0, n = 1, \beta^* = 0,25, \gamma = 7,95, Le = 10$. — stable steady states, – – – unstable steady states, • • • stable limit cycles, – • – unstable limit cycles.

Particular characteristics of the homoclinic orbit have now to be discussed. This orbit is formed as a result of the stable limit cycle joining both stable and unstable manifolds of the saddle, and surrounds an unstable stationary point. Therefore, the homoclinic orbit has to be internally stable; otherwise, the area inside the orbit would not possess any attractor, which is impossible on both mathematical and physical grounds. The nature of the external side of the homoclinic orbit was analysed using the phase diagram (Fig. 5).

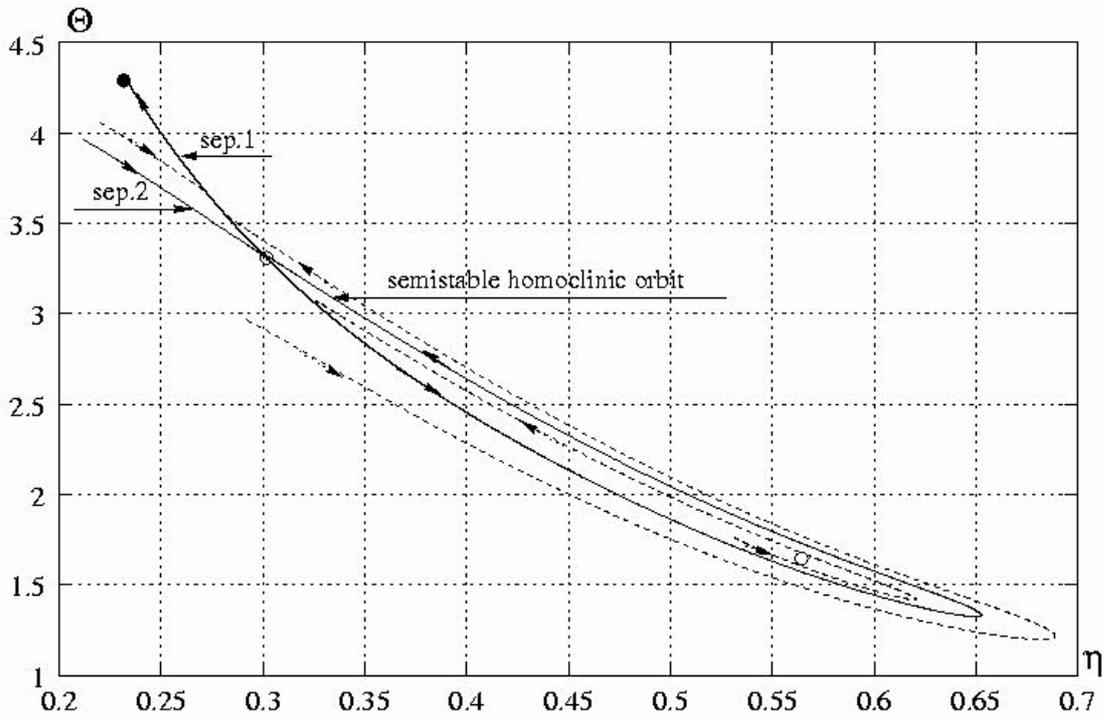

Fig.5.    Phase diagram: $a = 0, n = 1, \beta^* = 0,25, \gamma = 7,95, Le = 10, \theta_0 = 0,548621$. • stable steady state, ∘ unstable steady state, $--\rightarrow$ trajectories.

It is obvious from this diagram that this orbit is externally unstable. Consequently, the homoclinic orbit formed is semistable, i.e., it is stable on the inside and unstable on the outside (Fig. 5). It can be seen from the phase diagram that the trajectories originating inside the homoclinic orbit, when moving closer to the orbit finally meet the orbit and thus tend towards the saddle point (obviously, for an infinitely long period of time). This is quite a surprising conclusion, as a region has thus been created with a special property: trajectories originating from this region can reach an unstable saddle point (or, strictly speaking, can approach this point at an infinitely short distance). This phenomenon would be absolutely impossible without the existence of the above–mentioned homoclinic orbit, as otherwise the reaching of the saddle would be feasible only along the stable manifolds (which form the lines).

The dynamic system formed has therefore two attractors, one attractor proper (i.e. the stable stationary point on the upper branch of the solution diagram) and the other attractor which is reached indirectly via the homoclinic orbit (i.e. the saddle point). The latter attractor can be regarded as such only for trajectories originating from the inside of the region bounded by the homoclinic orbit.

For the Thiele moduli $\theta_0 = 0.548621 - 0.5489$ (Fig. 4) a region exists characterized by two unstable steady states (a saddle and the lower steady state) and one stable steady state (the upper steady state). This is an interesting example of the existence of a region of multiple steady states with only one stable steady state.

It follows from the bifurcation diagram (Fig. 2) that with a further increase in $\gamma$ above $\gamma = 7.957$ the loci of the Hopf bifurcation enter the region of multiple solutions, crossing the line of the limit points. Starting from the point at which the line of the Hopf bifurcations crosses that of the homoclinic bifurcations inside the region of multiple solutions, a basic change in the shape of the solution diagram can be seen (Fig. 6).

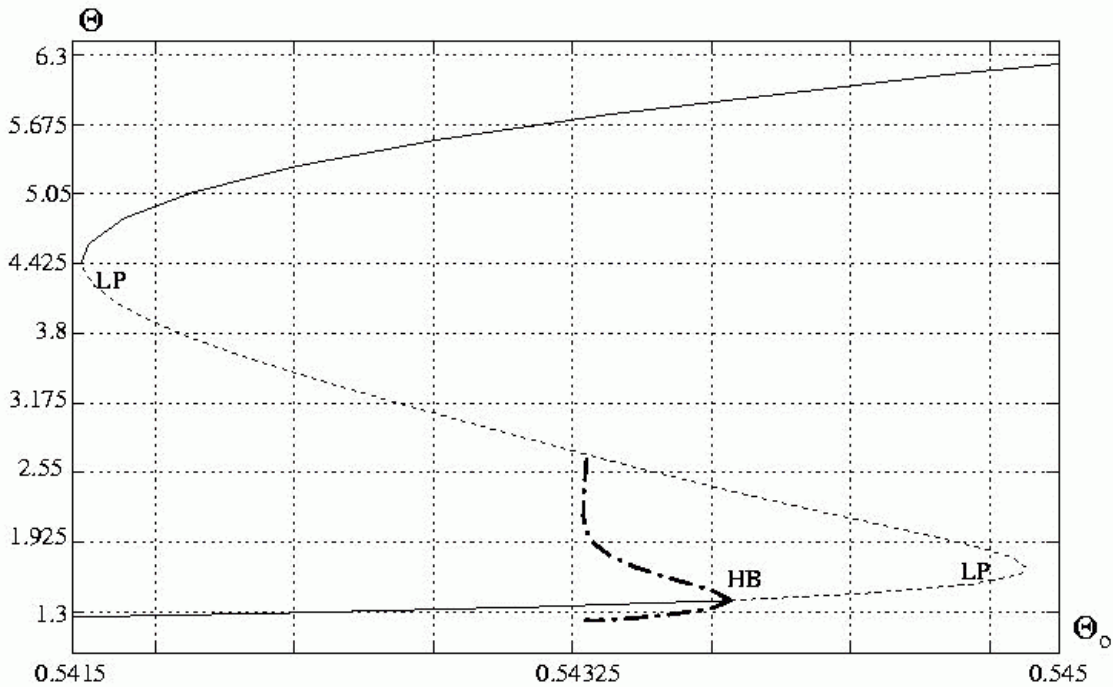

Fig.6.    Solution   diagram:   $a = 0, n = 1, \beta^* = 0,25, \gamma = 8,0, Le = 10$.   —   stable   steady   states, – – – unstable steady states, – • – unstable limit cycle.

Due to the small distance between the point of the Hopf bifurcation and the line of saddle points, the branch of periodic solutions is unable to develop fully and is lacking a part of the branch with stable limit cycles. Therefore, upon meeting the saddle point it forms an unstable homoclinic orbit illustrated in the phase diagram (Fig. 7).

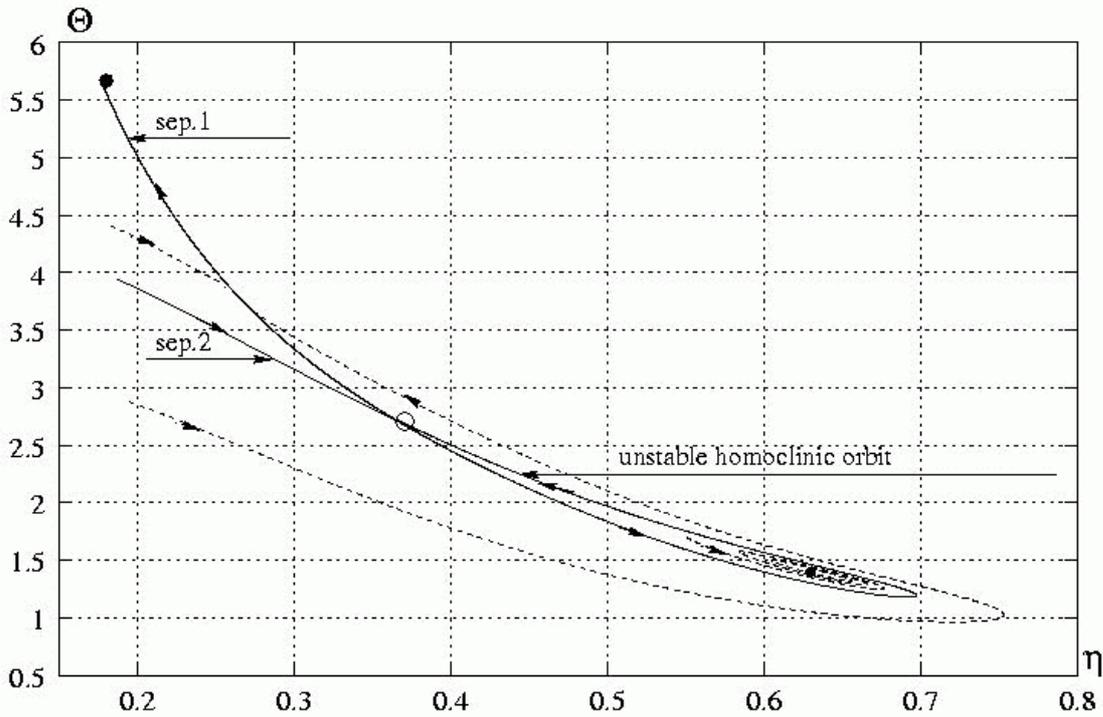

Fig.7.     Phase diagram: $a = 0, n = 1, \beta^* = 0,25, \gamma = 8,0, Le = 10, \theta_0 = 0,54332$. • stable steady state, ○
unstable steady state, $- - \rightarrow$ trajectories.

Consequently, the dynamic system shown in this diagram has two stable steady states (on both the upper and the lower branch of the steady states), along with a saddle point. The unstable homoclinic orbit plays here the role of a closed separatrix which divides the whole region into two subregions with different attractors (inside and outside the homoclinic orbit).

A further increase in the Thiele modulus leads to the break-up of the homoclinic orbit and the formation of an unstable limit cycle surrounding a stable steady state (Fig. 6). This cycle gradually diminishes and finally meets the Hopf bifurcation point $(\theta_0 = 0.5437)$, thus destabilizing this steady state. For $\theta_0 > 0.5437$ two unstable steady states exist in the region of multiple solutions (a saddle and the lower steady state) along with a single stable state (the upper steady state), similarly as in Fig. 4.

With an increase in $\gamma$ the Hopf bifurcation point, together with the homoclinic bifurcation point, moves towards the limit point (Fig. 1). When the meeting with the limit point (LP) occurs (this takes place at $\gamma = 8.55$), the Hopf bifurcation point disappears due to the so-called "double-zero singularity".

For Le = 50 the system (whose dynamics is illustrated in Figs 8-15) becomes by far more complex than that analysed previously, and reveals some additional interesting bifurcation phenomena.

It follows from the analysis of the bifurcation diagrams (Figs 8 and 9) that for $\gamma$ less than 7.95 a region of single steady states exists in which two lines of the Hopf bifurcation points occur. These lines define the ranges of the Thiele moduli associated with the feasible periodic solutions, similarly as in the previous case.

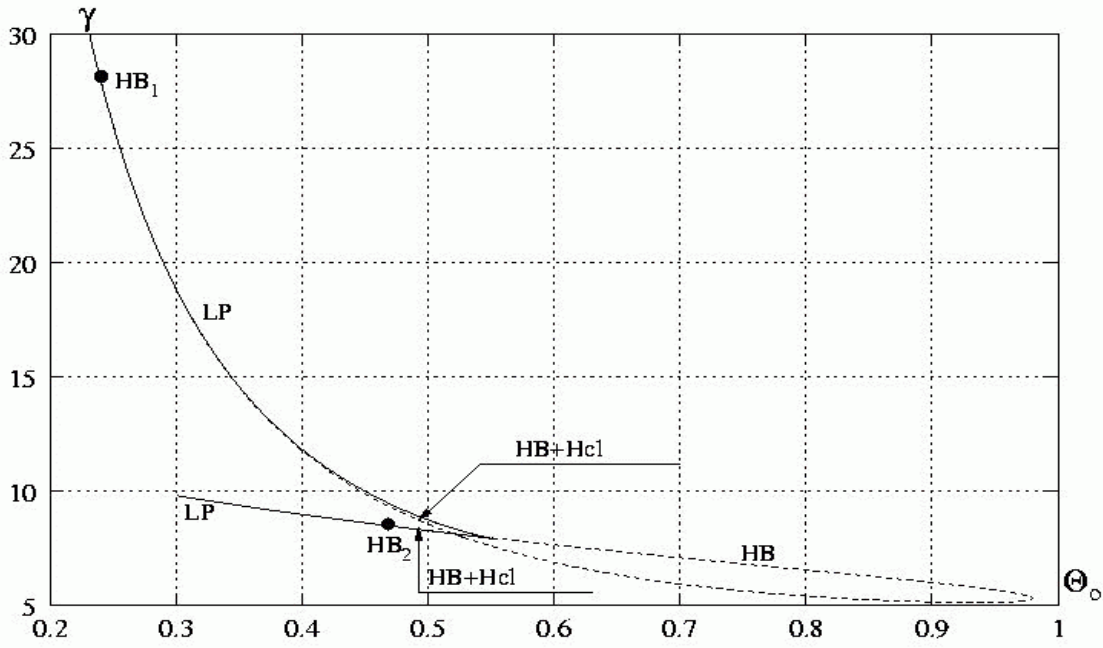

Fig.8. Bifurcation diagram: $a = 0, n = 1, \beta^* = 0,25, Le = 50$. — LP – locus of limit points, – – – HB – locus of Hopf bifurcation points, $\cdots$ Hcl – locus of homoclinic bifurcation points.

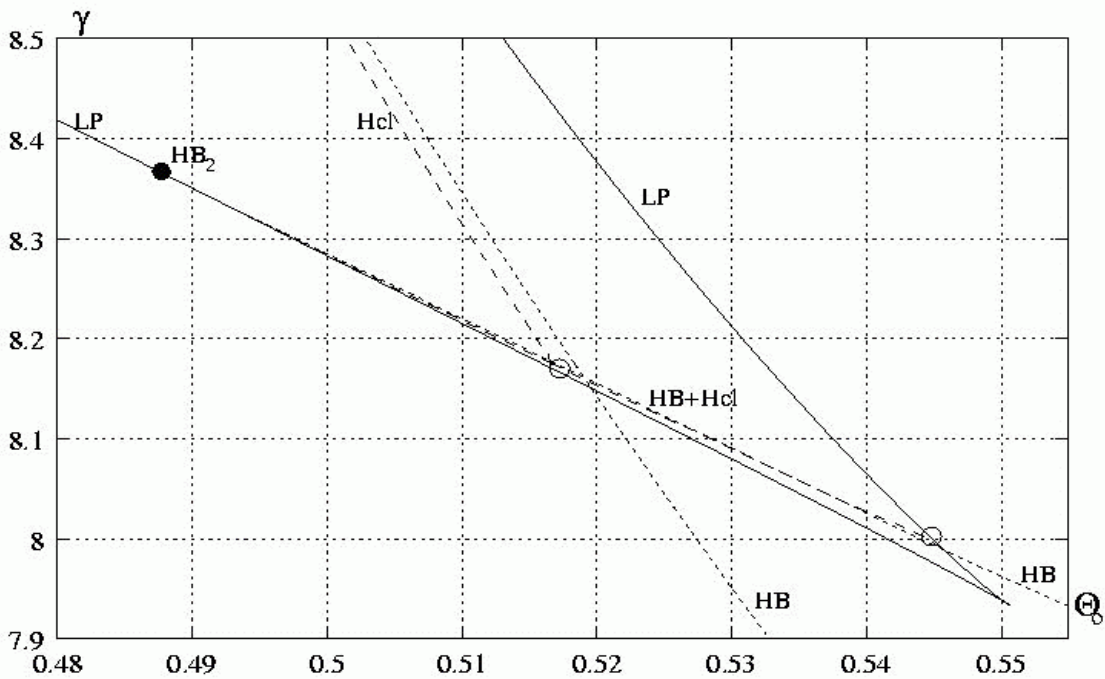

Fig.9. Bifurcation diagram: enlargement of a part of Fig. 8.

For $\gamma$ equal to about 7.94 a region of multiple steady states is formed entirely surrounded by a stable limit cycle. Such a system is illustrated in the solution diagram (Fig. 10) valid for $\gamma = 7.94 - 7.99$.

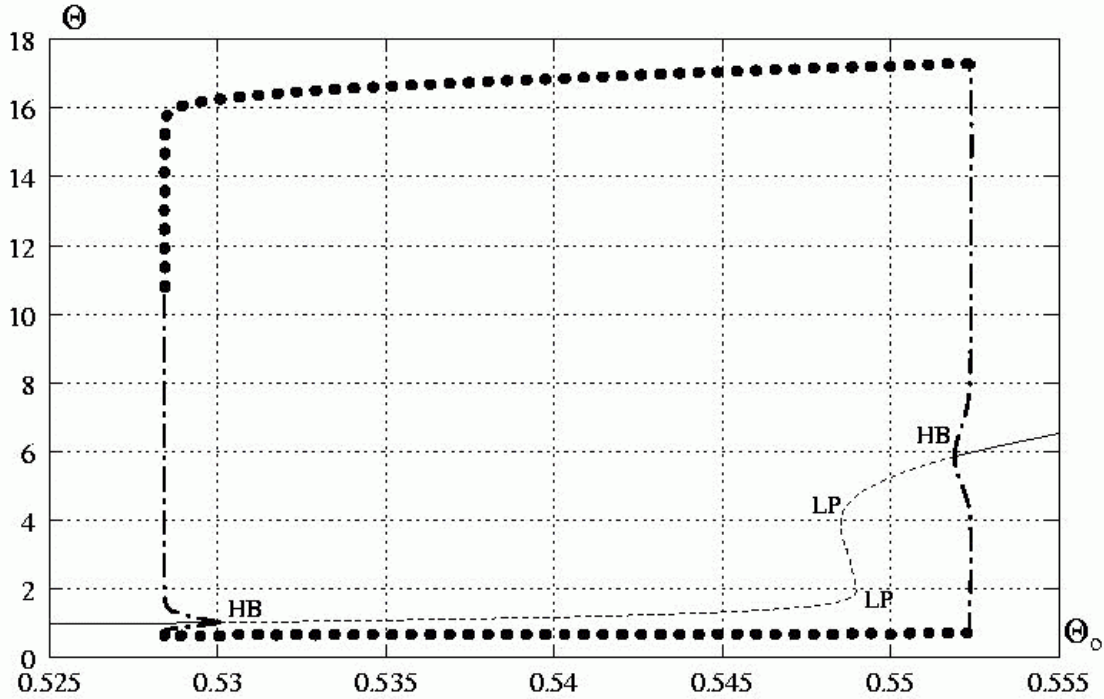

Fig.10. Solution diagram: $a = 0, n = 1, \beta^* = 0,25, \gamma = 7,95, \mathrm{Le} = 50$. — stable steady states, – – – unstable steady states, • • • stable limit cycles, – • – unstable limit cycles.

With an increase in $\gamma$ the Hopf bifurcation point, located on the upper branch of steady states, moves towards the upper limit point, while the Hopf bifurcation point located on the lower branch of steady states moves towards the lower limit point. An important moment is that when the limit cycle meets the lower limit point (due to the above-mentioned shift of the upper point of the Hopf bifurcation), and then when the Hopf bifurcation point enters the region of multiple solutions (which occurs at $\gamma \cong 7.99$).

At this point, as can be seen from the bifurcation diagram (Fig. 9), a homoclinic bifurcation line is initiated which, with increasing $\gamma$, is contained within the region of multiple steady states.

The solution diagram typical of the range of $\gamma = 7.99 - 8.16$ is shown in Fig. 11.

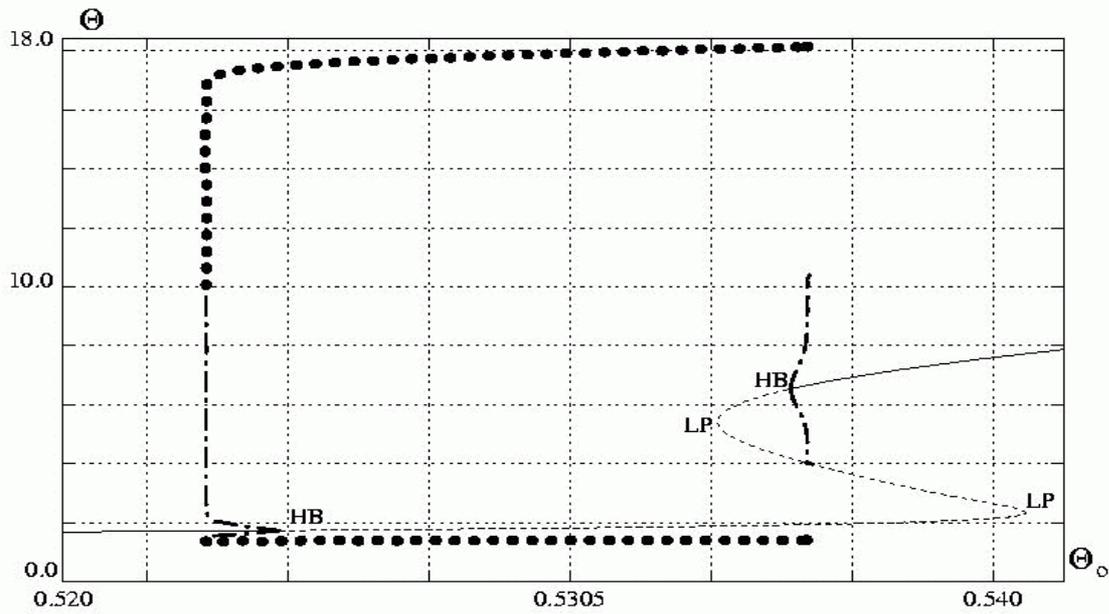

Fig.11. Solution diagram: $a = 0, n = 1, \beta^* = 0,25, \gamma = 8,05, Le = 50$. —— stable steady states, – – – unstable steady states, • • • stable limit cycles, – • – unstable limit cycles.

A more detailed illustration of the phenomenon of homoclinic bifurcation over the above range of the parameter $\gamma$ is given in the phase diagram for $\theta_0 = 0.53632$ (Fig. 12).

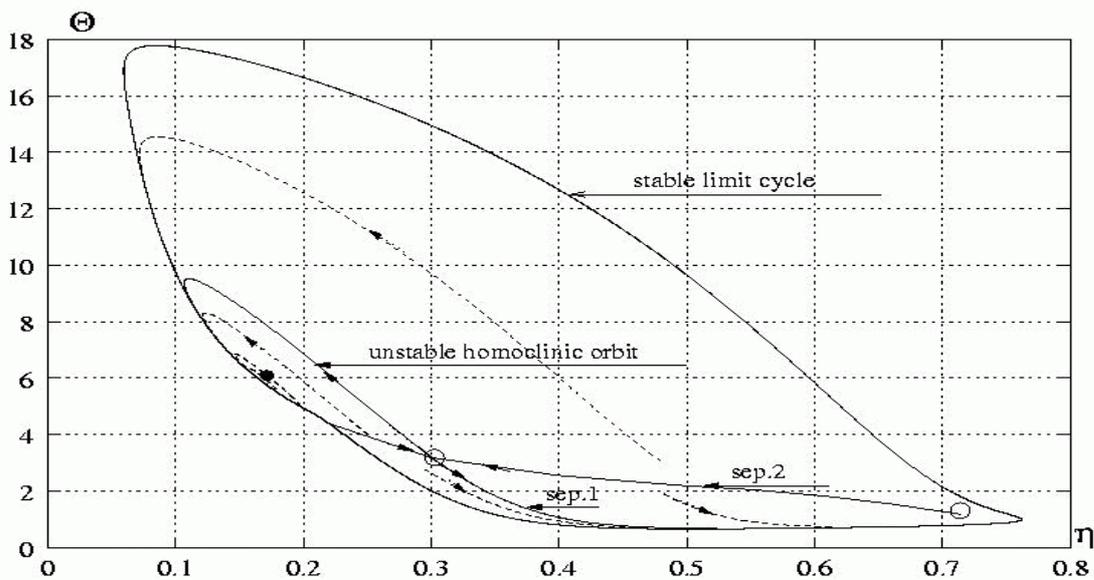

Fig.12. Phase diagram: $a = 0, n = 1, \beta^* = 0,25, \gamma = 8,05, Le = 50, \theta_0 = 0,53632$. • stable steady state, ∘ unstable steady state, – – → trajectories.

The analysis of this diagram reveals the existence of a large stable limit cycle, inside which an unstable homoclinic orbit is located that surrounds a stable stationary point. Outside the homoclinic orbit (obviously, over the region of the large stable limit cycle) an unstable stationary point is situated. The system discussed has therefore two attractors: a stable limit cycle and a stable stationary point inside the homoclinic orbit. The homoclinic orbit thus divides the region inside the stable limit cycle into two subregions. Over a certain distance the stable limit cycle and the unstable homoclinic orbit run so close to one another that, based on the phase diagram (Fig. 12), one might have the impression that they have merged altogether which, obviously, does not take place. Similarly, the visible trajectory clearly tends towards the stable stationary point; this is fully corroborated by numerical calculations.

Around the value of $\gamma \cong 8.16$ (cf. Fig. 9) the Hopf bifurcation point moving along the lower branch of steady states enters the region of multiple solutions, thus leading to the formation of a second line of the points of homoclinic bifurcation inside this region. Consequently, in the solution diagram (Fig. 13) two homoclinic orbits form for the given values of the Thiele modulus.

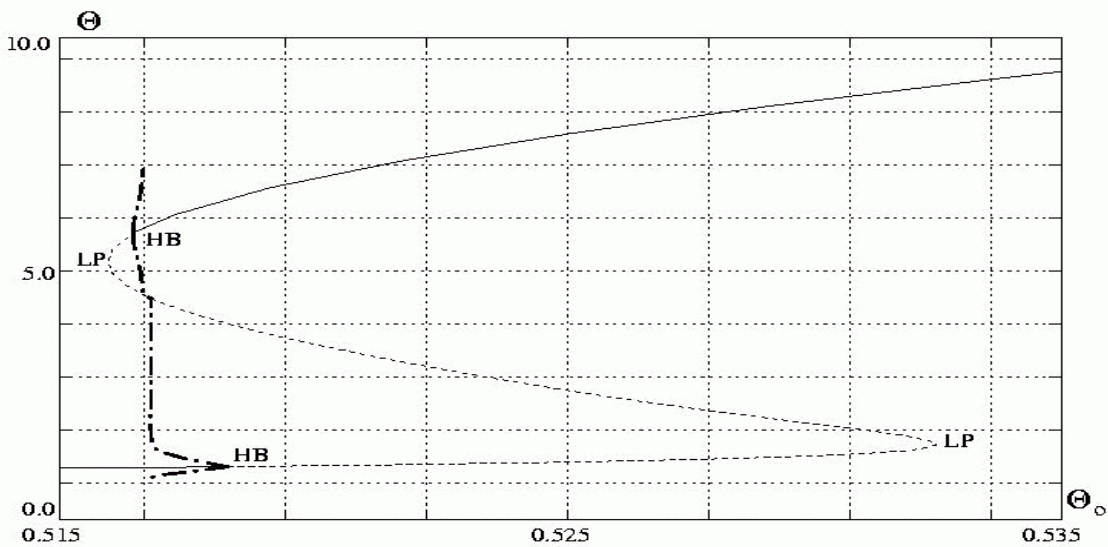

Fig.13. Solution diagram: $a = 0, n = 1, \beta^* = 0,25, \gamma = 8,175, Le = 50$. —— stable steady states, – – – unstable steady states, – • – unstable limit cycles.

A further increase in $\gamma$ (Fig. 9) leads to the merger of the Hopf bifurcation point located on the upper branch of steady states with the upper limit point $\left(\gamma \cong 8.37\right)$, to the disappearance of this point of the Hopf bifurcation and, finally, to the disappearance of periodic solutions.

On the other hand, the Hopf bifurcation point moving along the lower branch of steady states towards the lower limit point generates branches of periodic solutions of the form illustrated in Fig. 14.

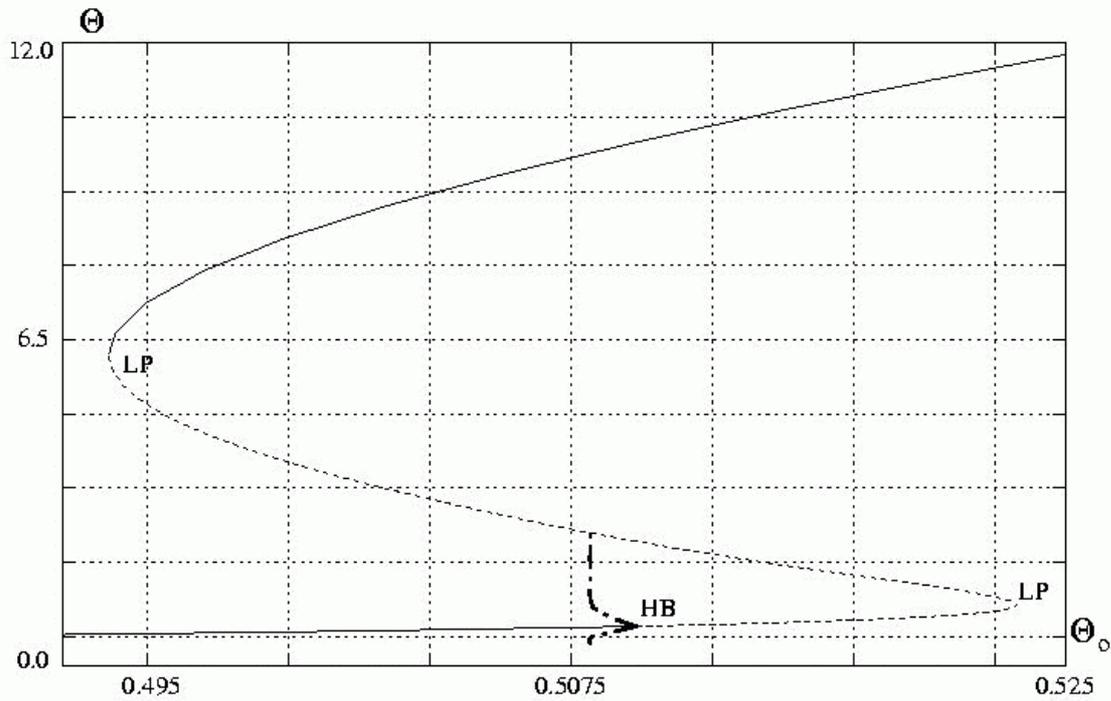

Fig.14.   Solution   diagram:   $a = 0, n = 1, \beta^* = 0{,}25, \gamma = 8{,}32, \text{Le} = 50$ .   —   stable   steady   states, – – – unstable steady states, – • – unstable limit cycles.

The phase diagram (Fig. 15) corresponding to this solution diagram shows the system dynamics for $\theta_0 = 0.50972$, i.e., for such a value of the Thiele modulus for which homoclinic bifurcation takes place.

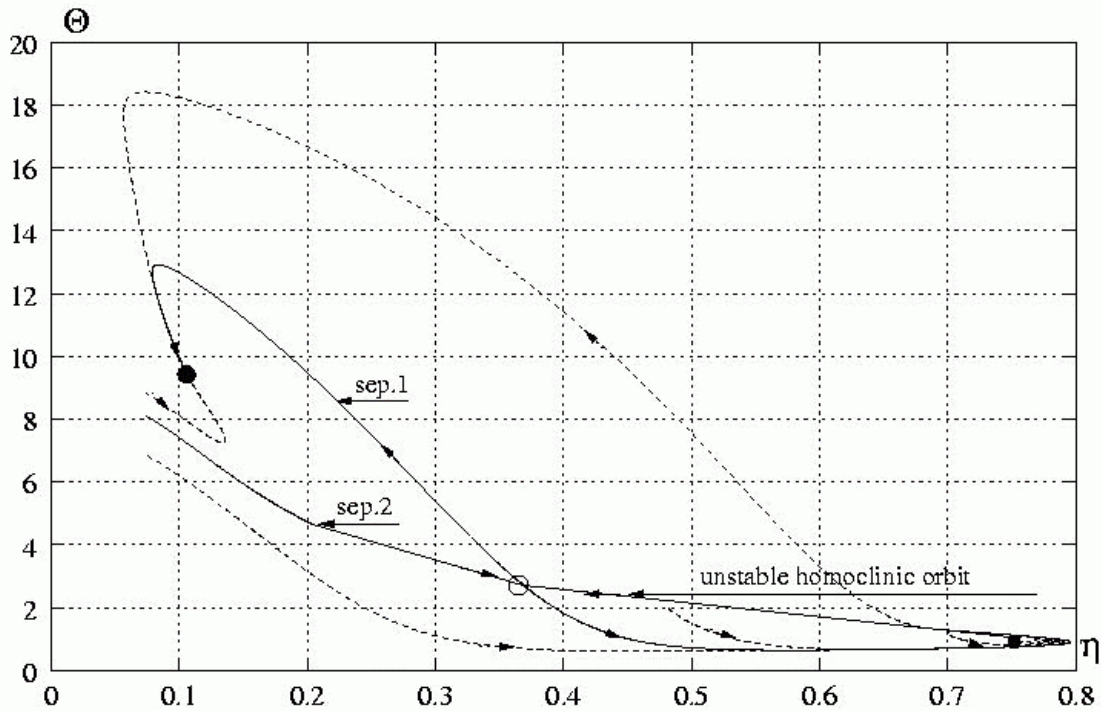

Fig.15.    Phase diagram: $a = 0, n = 1, \beta^* = 0{,}25, \gamma = 8{,}32, Le = 50, \theta_0 = 0{,}50972$. • stable steady state, ○ unstable steady state, $-- \rightarrow$ trajectories.

In this diagram, apart from two stable steady states an unstable homoclinic orbit can be seen. The trajectory shown in this diagram passes outside the unstable homoclinic orbit and approaches this orbit so closely that it seems to glide along the orbit; then, after covering a certain distance, it leaves the orbit and tends towards the stable stationary point.

It follows from the analysis of the solution diagrams (Figs 11, 13 and 14) and the corresponding phase diagrams (Figs 12 and 15) that the homoclinic orbits for the system studied are unstable. This is due to the fact that the branches of periodic solutions, which correspond to the subcritical Hopf bifurcation include only a part of the branch containing unstable limit cycles as, before the part of the branch with the stable limit cycles becomes fully developed, it meets the saddle point to form a homoclinic orbit.

For $\gamma \cong 28$ the Hopf bifurcation point located on the lower branch of steady states joins the limit point and vanishes altogether.

## 4. CONCLUSIONS

The analysis performed and the simulations carried out in the present study reveal the possibility of the occurrence of the homoclinic bifurcation which forms a homoclinic orbit during the phenomena taking place in a catalyst pellet.

A method is developed which makes it possible to determine, in the bifurcation diagram, the loci of the homoclinic bifurcation (Hcl), together with the lines of the limit points (LP) and the loci of the Hopf bifurcation. It is clearly seen that the formation of the homoclinic bifurcation depends upon the coexistence of periodic solutions and the region of multiple steady states and, in particular, upon the existence of the unstable branch of steady states with all the characteristics of saddle points in the neighbourhood of the branch of periodic solutions. It is only under these conditions that the limit cycle can join the saddle point to form the homoclinic orbit. Furthermore, the bifurcation diagrams define such a set of parameters for which homoclinic bifurcations can appear.

It is found that two types of homoclinic orbits can form, namely, an unstable orbit and a semistable orbit. The unstable homoclinic orbit appears as a result of the merger between an unstable limit cycle and a saddle point and surrounds the stable stationary point. This orbit plays therefore the role of a closed separatrix, dividing the region of state variables in the phase diagram into two subregions with two attractors. Two stable stationary points or a stable stationary point inside the orbit and a stable limit cycle outside the homoclinic orbit can be these attractors.

On the other hand, the merger of a stable limit cycle and a saddle point leads to the formation of a semistable homoclinic orbit which is stable on the inside and unstable on the outside. The internal stability of the orbit is due to the fact that the orbit surrounds an unstable stationary point. The absence of stability of the internal side of the homoclinic orbit would mean that the region inside the orbit does not have any attractor, which is impossible form the standpoint of the principles of dynamics.

A further highly interesting dynamic phenomenon is associated with the existence of a semistable homoclinic orbit. The internal stability of this orbit leads to the fact that the trajectories originating from the region surrounded by the semistable homoclinic orbit tend to this orbit and, consequently, to the unstable saddle point (obviously, after an infinitely long time). A region of state variables is thus created with a special property: the trajectories issuing from this region can reach the unstable stationary point – a phenomenon absolutely impossible without the existence of the homoclinic orbit.

It has also to be stressed that all the points of the Hopf bifurcation determined in the present study are the points of subcritical bifurcation. This leads to the possibility of the occurrence, under strictly defined conditions, of unstable homoclinic orbits.

## NOTATION

a – geometric factor defining the shape of the pellet (a = 0, slab; a = 1, cylinder; a = 2, sphere)

$Bi_H = \dfrac{\alpha L}{\lambda_e}$        – Biot number for heat transfer

$Bi_M = \dfrac{kL}{D_e}$        – Biot number for mass transport

$C_A$        – concentration of species A, $kmol \cdot m^{-3}$

$D_e$        – effective diffusivity, $m^2 \cdot s^{-1}$

E        – activation energy, $kJ \cdot kmol^{-1}$

HB        – Hopf bifurcation point

Hcl        – homoclinic bifurcation point

$\Delta H$        – heat of reaction, $kJ \cdot kmol^{-1}$

k        – mass transfer coefficient, $m \cdot s^{-1}$

L        – characteristic dimension of the pellet, m

$Le = \dfrac{\alpha L(a+1)}{(\rho c)_k D_e}$        – modified Lewis number

$Le' = \dfrac{\lambda_e}{(\rho c)_k D_e}$        – Lewis number

LP        – limit point

| | |
|---|---|
| n | – reaction order |
| $r(C_A, T)$ | – reaction rate, $kmol \cdot m^{-3} \cdot s^{-1}$ |
| $R(y,z) = \dfrac{r_A(C_A, T)}{r_A(C_{A0}, T_0)}$ | – dimensionless rate of reaction |
| T | – temperature, K |
| t | – time, s |
| V | – volume of the pellet, $m^3$ |
| $x = \dfrac{\rho}{\rho_0}$ | – dimensionless coordinate in the pellet |
| $y = \dfrac{C_A}{C_{A0}}$ | – dimensionless concentration |
| $z = \dfrac{T}{T_0}$ | – dimensionless temperature |

## GREEK LETTERS

| | |
|---|---|
| $\alpha$ | – heat transfer coefficient, $W \cdot m^{-2} \cdot K^{-1}$ |
| $\beta = \dfrac{(-\Delta H) D_e C_{A0}}{T_0 \lambda_e}$ | – dimensionless Prater number |
| $\beta^* = \dfrac{2\beta(a+1)}{Bi_H(n+1)}$ | – dimensionless parameter |
| $\gamma = \dfrac{E}{R_g T_0}$ | – dimensionless activation energy |
| $\theta_0 = \dfrac{\phi_0}{a+1}\left(\dfrac{n+1}{2}\right)^{1/2}$ | – modified Thiele modulus |
| $\theta^2 = \theta_0^2 \exp\left[\gamma\left(1 - \dfrac{1}{z}\right)\right]$ | – transformed dimensionless temperature |
| $\lambda_e$ | – effective conductivity, $W \cdot m^{-1} \cdot K^{-1}$ |
| $\nu$ | – stoichiometric coefficient $(\nu_A = -1)$ |
| $(\rho c)_k$ | – heat capacity of catalyst, $kJ \cdot m^{-3} \cdot K^{-1}$ |
| $\tau = \dfrac{D_e t}{L^2}$ | – dimensionless time |
| $\phi_0^2 = \dfrac{L^2 r(C_{A0}, T_0)}{D_e C_{A0}}$ | – Thiele modulus |

## SUBSCRIPTS

| | |
|---|---|
| 0 | – refers to conditions in bulk phase |
| z | – refers to surface conditions |